\documentstyle[twoside,11pt,leqno]{article}

\def\ind{\textrm{ind}} \def\X{{\cal X}}  \def\H{{\cal H}} 
\def\d{\rm{des}}
\def\i{\rm{iso}} \def\a{\rm{asc}}
\def\n{{\texttt N}} \def\N{{\texttt{N}}}
\def\C{\texttt{C}} \def\iso{\textrm{iso}}
\def\p{{(\cal P)}}
\def\P{{(\cal {HP})}}

\def\B{B({\cal H})} \def\b{B({\cal X})}
\def\asc{ \textrm{asc}} \def\dsc{ \textrm{dsc}}

\newtheorem{df}{Definition}[section]
\newtheorem{thm}[df]{Theorem} \newtheorem{pro}[df]{Proposition}
 
\newtheorem{rema}[df] {Remark} 
\def\sfstp{{\hskip-1em}{\bf.}{\hskip1em}}

\def\subject#1{\renewcommand{\thefootnote}{}\footnote
{AMS(MOS) subject classification (2010). Primary: {#1}}}

\def\keywords#1{\renewcommand{\thefootnote}{}\footnote
{Keywords: {#1}}}

\def\enddemo{\qed \endtrivlist} \expandafter\let\csname
enddemo*\endcsname=\enddemo

\def\qedsymbol{\ifmmode\bgroup\else$\bgroup\aftergroup$\fi
\vcenter{\hrule\hbox{\vrule
height.5em\kern.5em\vrule}\hrule}\egroup}
\def\qed{\ifmmode\else\unskip\nobreak\fi\quad\qedsymbol}

\title{\bf  Compact perturbations and consequent hereditarily polaroid operators}
\author{\normalsize B.P. Duggal }
\date{}

\begin{document}

\maketitle \thispagestyle{empty} \vskip-16pt

\subject{47A10, 47A55, 47A53, 47B40} \keywords{Banach space
operator, Hilbert space, compact perturbation, SVEP, hereditarily
polaroid, Fredholm }

\begin{abstract}A Banach space operator $A\in\b$ is polaroid,
$A\in\p$, if the isolated points of the spectrum $\sigma(A)$ are
poles of the operator; $A$ is hereditarily polaroid, $A\in\P$, if
every restriction of $A$ to a closed invariant subspace is polaroid.
Operators $A\in\P$ have SVEP on $\Phi_{sf}(A)=\{\lambda: A-\lambda$
is semi Fredholm $\}$: This, in answer to a question posed by Li and
Zhou \cite[Problem 5.5]{LZ}, proves the necessity of the condition
$\Phi_{sf}^+(A)=\emptyset$. A sufficient condition for $A\in\b$ to
have SVEP on $\Phi_{sf}(A)$ is that its component
$\Omega_a(A)=\{\lambda\in\Phi_{sf}(A): \ind(A-\lambda)\leq 0\}$ is
connected. We prove: If $A\in\B$ is a Hilbert space operator, then a
necessary and sufficient condition
 for there to exist a compact operator $K\in\B$ such that $A+K\in\P$ is that
 $\Omega_a(A)$ is connected.
\end{abstract}

%%%%%%%%%%%%%%%%%%%%%%%%%%%%%%%%%%%%%%%% SECTION 1

\section {\sfstp Introduction} Let $\b$ (resp., $\B$) denote the
algebra of operators, equivalently bounded linear transformations,
on a complex infinite dimensional Banach (resp., Hilbert) space into
itself. For an operator $A\in\b$, let $\i\sigma(A)$ denote the
isolated points of the spectrum $\sigma(A)$, let $\asc(A)$ (resp.,
$\dsc(A)$) denote the ascent (resp., descent) of $A$ and let
$A-\lambda$ denote $A-\lambda I$. A point $\lambda\in\i\sigma(A)$ is
a pole (of the resolvent) of $A$, equivalently $A$ is polar at
$\lambda$, if $\asc(A-\lambda)=\dsc(A-\lambda)<\infty$. The operator
$A$ is {\em polaroid} if it is polar at every
$\lambda\in\i\sigma(A)$, and it is {\em hereditarily polaroid} if
every restriction  $A|_M$ of $A$ to an (always closed) invariant
subspace $M$ of $A$ is polaroid. Polaroid operators, and their
perturbation by commuting compact perturbations, have been studied
by a number of authors in the recent past (see \cite{{AA}, {A1},
{D}, {D1}, {LZ}} for a sample). For example, if $N\in\b$ is a
nilpotent operator which commutes with $A\in\b$, then $A$ is
polaroid if and only if $A+N$ is polaroid \cite[Theorem 2.6(b)]{D1}.
This however does not extend to non-nilpotent quasinilpotent
commuting operators: Consider for example the trivial operator
$A=0\in\b$ and a non-nilpotent quasinilpotent $Q\in\b$. The
perturbation of a polaroid operator by a compact operator may or may
not effect the polaroid property of the operator. For example, if
$U\in\B$ is the forward unilateral shift, $A=U\oplus U^*$ and $K$ is
the compact operator $K= \left(\begin{array}{clcr} 0 & 1-UU^*\\0 &
0\end{array}\right)$, then both $A$ and $A+K$ are polaroid (for the
reason that $\i\sigma(A)=\emptyset$ and $A+K$ is a unitary);
trivially the identity operator $1$ is polaroid, but its
perturbation $1+Q$ by a compact quasinilpotent operator is not
polaroid.

\

An interesting problem, recently considered by Li and Zhou
\cite{LZ}, is the following: Given an operator $A\in\B$, do there
exist compact operators $K_0, K\in\B$ such that (i) $A+K_0$ is
polaroid and (ii) $A+K$ is not polaroid. The answer to both these
problems is an emphatic ``yes" (see \cite[Theorems 1.4 and 1.5]{LZ}.
The argument used to prove these results ties up with the work of
Herrero and his co-authors \cite{{He}, {He1}, {HTW}, {AFHV}}, Ji
\cite{Ji}, and Zhu and Li \cite{ZL}. A natural extension of this
problem is the question of whether there exist compact operators
$K_0, K\in\B$ such that (i)' $A+K_0$ is hereditarily polaroid and
(ii)' $A+K$ is not hereditarily polaroid. Here the answer to (ii)'
is a ``yes" \cite[Theorem 5.2]{LZ}, but there is caveat to the
answer to (i)' - the answer is ``yes if the set
$\Phi_{sf}^+(A)=\{\lambda\in\sigma(A): A-\lambda$ is semi- Fredholm
and $\ind(A-\lambda)>0\}=\emptyset$". The authors of \cite{LZ} leave
the problem of a straight ``yes or no" answer to (i)' open. This
note considers this problem to prove that if $A, K\in\b$ with $K$
compact, then $A+K$ hereditarily polaroid implies
$\Phi_{sf}^+(A)=\emptyset$. Indeed, we prove that if $A\in\B$, then
there exists a compact $K\in\B$ such that $A+K$ is hereditarily
polaroid if and only if $A$ has SVEP, {\em the single-value
extension property}, on $\Phi_{sf}(A)$. A sufficient condition for
operators $A\in\b$ to have SVEP on $\Phi_{sf}(A)$ is that the
component $\Omega_a(A)=\{\lambda\in\Phi_{sf}(A): \ind(A-\lambda)\leq
0\}$ is connected. We prove that for an operator $A\in\B$, a
necessary and sufficient condition
 for there to exist a compact operator $K\in\B$ such that $A+K\in\P$ is that
 $\Omega_a(A)$ is connected.

\section {\sfstp  Complementary results}  We start by introducing our notation and terminology.
We shall denote the class of polaroid operators by $\p$ and the
subclass of hereditarily polaroid operators by $\P$. The boundary of
a subset $S$ of the set $\C$ of complex numbers will be denoted by
$\partial{S}$. An operator $A\in\b$ has SVEP, {\em the single-valued
extension property}, at a point $\lambda_0\in\C$ if for every open
disc ${\mathcal D}_{\lambda_0}$ centered at $\lambda_0$ the only
analytic function $f:{\mathcal D}_{\lambda_0}\longrightarrow \X$
satisfying $(A-\lambda)f(\lambda)=0$ is the function $f\equiv 0$.
(Here, as before, we have shortened $A-\lambda 1$ to $A-\lambda$.)
Evidently, every $A$ has SVEP at points in the resolvent
$\rho(A)=\C\setminus \sigma(A)$ and the boundary $\partial\sigma(A)$
of the spectrum $\sigma(A)$. {\em We say that $T$ has SVEP on a set
$S$ if it has SVEP at every $\lambda\in S$.} The {\em ascent of
$A$}, $\asc(A)$ (resp. {\em descent of $A$}, $\dsc(A)$), is the
least non-negative integer $n$ such that $A^{-n}(0)=A^{-(n+1)}(0)$
(resp., $A^n(\X)=A^{n+1}(\X)$): If no such integer exists, then
$\asc(A)$, resp. $\dsc(A)$, $=\infty$. It is well known  that
$\asc(A)<\infty$ implies $A$ has SVEP at $0$, $\dsc(A)<\infty$
implies $A^*$ ($=$ the dual operator) has SVEP at $0$, finite ascent
and descent for an operator implies their equality, and that a point
$\lambda\in\sigma(A)$ is a pole (of the resolvent) of $A$ if and
only if $\asc(A-\lambda)=\dsc(A-\lambda)<\infty$ (see \cite{{A},
{H}, {LN}}).

\

 An
operator $A\in \b$ is: {\em upper semi--Fredholm at $\lambda\in\C$},
$\lambda\in\Phi_{uf}(A)$ or $A-\lambda\in\Phi_{uf}(\X)$, if
$(A-\lambda)(\X)$ is closed and the deficiency
 index
 $\alpha(A-\lambda)=\rm{dim}((A-\lambda)^{-1}(0))<\infty$; {\em
lower semi--Fredholm at $\lambda\in\C$}, $\lambda\in\Phi_{lf}(A)$ or
$A-\lambda\in\Phi_{lf}(\X)$, if $\beta(A-\lambda)=\rm{dim}(\X/
(A-\lambda)(\X))<\infty$.
 $A$ is semi--Fredholm, $\lambda\in\Phi_{sf}(A)$ or  $A-\lambda\in\Phi_{sf}(\X)$,
 if $A-\lambda$ is either upper or lower semi--Fredholm,  and $A$ is Fredholm,
 $\lambda\in\Phi(A)$ or $A-\lambda\in\Phi(\X)$, if $A-\lambda$ is both upper and
 lower semi--Fredholm. The index of a semi--Fredholm operator  is
the integer, finite or infinite, $\ind(A)=\alpha(A)-\beta(A)$.
Corresponding to these classes of one sided Fredholm operators, we
have the following spectra: The {\em upper Fredholm spectrum}
$\sigma_{uf}(A)$ of $A$ defined by
$\sigma_{uf}(A)=\{\lambda\in\sigma(A):A-\lambda\notin\Phi_{uf}(\X)\}$,
and the {\em lower Fredholm spectrum} $\sigma_{lf}(A)$ of $A$
defined by
$\sigma_{lf}(A)=\{\lambda\in\sigma(A):A-\lambda\notin\Phi_{lf}(\X)\}$.
The {\em Fredholm spectrum} $\sigma_f(A)$ of $A$ is the set
$\sigma_f(A)=\sigma_{uf}(A)\cup\sigma_{lf}(A)$, and the {\em Wolf
spectrum} $\sigma_{ulf}(A)$ of $A$ is the set
$\sigma_{ulf}(A)=\sigma_{uf}(A)\cap\sigma_{lf}(A)$. $A\in\b$ is
 Weyl (at $0$) if it is  Fredholm with $\ind(A)=0 $. It is well known
that a semi- Fredholm operator $A$ (resp., its conjugate operator
$A^*$) has SVEP at a point $\lambda$ if and only if
$\asc(A-\lambda)<\infty$ (resp., $\dsc(A-\lambda)<\infty$)
\cite[Theorems 3.16, 3.17]{A}; furthermore, if $A-\lambda$ is Weyl
 , i.e., if $\lambda\in\Phi(A)$ and
$\ind(A-\lambda)=0$, then $A$ has SVEP at $\lambda$ implies
$\lambda\in\i\sigma(A)$ with
$\asc(A-\lambda)=\dsc(A-\lambda)<\infty$. The Weyl (resp., the upper
or approximate Weyl) spectrum of $A$ is the set \begin{eqnarray*} &
& \sigma_{w}(A)=\{\lambda\in\sigma(A):
\lambda\in\sigma_{f}(A)\hspace{2mm}\mbox{or}\hspace{2mm}\ind(A-\lambda)\neq
0\}\\ & & (\sigma_{aw}(A)=\{\lambda\in\sigma_a(A):
\lambda\in\sigma_{uf}(A)\hspace{2mm}\mbox{or}\hspace{2mm}\ind(A-\lambda)>
0\}).\end{eqnarray*} The Browder (resp., the upper or approximate
Browder) spectrum of $A$ is the set
\begin{eqnarray*}
& &
\sigma_{b}(A)=\{\lambda\in\sigma(A):\lambda\in\sigma_f(A)\hspace{2mm}
\mbox{or}\hspace{2mm} \a(A-\lambda)\neq \d(A-\lambda)\}\\ & &
(\sigma_{ab}(A)=\{\lambda\in\sigma_a(A):\lambda\in\sigma_{uf}(A)\hspace{2mm}
\mbox{or}\hspace{2mm} \a(A-\lambda)= \infty\}.
\end{eqnarray*} Clearly,
$\sigma_f(A)\subseteq\sigma_w(A)\subseteq\sigma_b(A)\subseteq\sigma(A)$
and
$\sigma_{uf}(A)\subseteq\sigma_{aw}(A)\subseteq\sigma_{ab}(A)\subseteq\sigma(A)$.

\

An operator $A\in\b$ is B-Fredholm (resp., upper B-Fredholm),
$A\in\Phi_{Bf}(\X)$ (resp., $\Phi_{uBf}(\X)$), if there exists an
integer $n\geq 1$ such that $A^n(\X)$ is closed and the induced
operator $A_{[n]}=A|_{A^n(\X)}$, $A_{[0]}=A$, is Fredholm (resp.,
upper semi Fredholm) in the usual sense. It is seen that if
$A_{[n]}\in\Phi_{sf}(\X)$ for an integer $n\geq 1$, then
$A_{[m]}\in\Phi_{sf}(\X)$ for all integers $m\geq n$: One may thus
define unambigiously the index of $A$ by
$\ind(A)=\alpha(A)-\beta(A)$ (see \cite{{BK}, {A1}, {AZ}}). The
B-Fredholm (resp., the upper B-Fredholm) spectrum of $A$ is the set
\begin{eqnarray*} & &\sigma_{Bf}(A)=\{\lambda\in\sigma(A):\lambda\notin\Phi_{Bf}(A)\}\\
& &
(\sigma_{uBf}(A)=\{\lambda\in\sigma(A):\lambda\notin\Phi_{uBf}(A)\})\end{eqnarray*}
and the B-Weyl (resp.,upper or approximate B-Weyl) spectrum of $A$
is the set \begin{eqnarray*} & &
\sigma_{Bw}(A)=\{\lambda\in\sigma(A):
\lambda\in\sigma_{Bf}(A)\hspace{2mm}\mbox{or}\hspace{2mm}\ind(A-\lambda)\neq
0\}\\ & & (\sigma_{uBw}(A)=\{\lambda\in\sigma(A):
\lambda\in\sigma_{uBf}(A)\hspace{2mm}\mbox{or}\hspace{2mm}\ind(A-\lambda)>
0\}).\end{eqnarray*} It is clear that
$\sigma_{Bw}(A)\subseteq\sigma_w(A)$ and
$\sigma_{uBw}(A)\subseteq\sigma_{aw}(A)$.

\

Let $H_0(A)$ and $K(A)$ denote, respectively, the {\em
quasinilpotent part}
$$H_0(A)=\{x\in\X:\lim_{n\rightarrow\infty}{||A^nx||^{\frac{1}{n}}=0}\}$$

 and the {\em analytic core} \begin{eqnarray*} & & K(A)=\{x\in\X:\hspace{1mm}\mbox{there
 exists a sequence}\hspace{2mm}
 \{x_n\}\subset\X\hspace{2mm}\mbox{and} \hspace{2mm}\delta>0
 \hspace{2mm}\mbox{for which}\\ & & x=x_0, Ax_{n+1}=x_n
 \hspace{2mm}\mbox{and}\hspace{2mm} ||x_n||\leq \delta^n ||x||
 \hspace{2mm}\mbox{for all}\hspace{2mm} n=0, 1, 2, \cdots
 \}\end{eqnarray*} of $A$. It is well known, \cite{A}, that
 $(A-\lambda)^{-p}(0)\subseteq H_0(A-\lambda)$, for all integers
 $p\geq 1$, and $(A-\lambda)K(A-\lambda)=K(A-\lambda)$  for all
 complex $\lambda$. A necessary and sufficient condition for
 $\lambda\in\i\sigma(A)$ to be a pole of $A$ is that
 $H_0(A-\lambda)=(A-\lambda)^{-p}(0)$ for some integer $p\geq 1$:
 This is seen as follows. If $\lambda\in\i\sigma(A)$, then (by the
 Riesz representation theorem \cite{{A}, {H}})\begin{eqnarray*} & &
 \X= H_0(A-\lambda)\oplus K(A-\lambda)= (A-\lambda)^{-p}(0)\oplus
 K(A-\lambda)\\ & \Longrightarrow & (A-\lambda)^p(\X)= 0\oplus
 (A-\lambda)^pK(A-\lambda)= K(A-\lambda)\\ & \Longrightarrow & \X=
 (A-\lambda)^{-p}(0)\oplus (A-\lambda)^p(\X)\\ & \Longrightarrow &
 \lambda \hspace{2mm}\mbox{is a pole of order}\hspace{2mm}
 p\hspace{2mm}\mbox{of}\hspace{2mm} A\\ & \Longrightarrow &
 H_0(A-\lambda)= (A-\lambda)^{-p}(0).\end{eqnarray*} For every
 $\lambda\notin\sigma_{Bw}(A)$ such that $A$ has SVEP at $\lambda$,
 $\asc(A-\lambda)<\infty$  (implying thereby that there
 exists an integer $p\geq 1$ such that
 $H_0(A-\lambda)=(A-\lambda)^{-p}(0)$) and $\lambda\in\i\sigma(A)$
 \cite{D3}. Hence $A$ has SVEP at $\lambda\notin\sigma_{Bw}(A)$
 implies $\lambda$ is a pole of $A$.

 If $\X=\H$ is a Hilbert space, and $A\in\B$ is such that $\lambda\in\Phi_{sf}(A)$,
 then {\em the minimal index of $A-\lambda$} is the integer $$
 \rm{min}\{\alpha(A-\lambda),
 \beta(A-\lambda)\}=\rm{min}\{\alpha(A-\lambda),
 \alpha(A-\lambda)^*\}.$$ It is well known that the function
 $\lambda\rightarrow \rm{min}.\ind(A-\lambda)$ is constant on every
 component of $\Phi_{sf}(A)$ (except perhaps for a denumerable subset
 without limit points in $\Phi_{sf}(A)$) \cite[Corollary 1.14]{He}

\section {\sfstp  Results} Given $A\in\b$, the {\em reduced minimum
modulus function} $\gamma(A)$ is the function $$\gamma(A) =
\rm{inf}_{x\notin A^{-1}(0)}\{\frac{||Ax||}{\rm{dist}(x,
A^{-1}(0))}\},$$ where $\gamma(A)=\infty$ if $A=0$. Recall that
$A(\X)$ is closed if and only if $\gamma(A)> 0$. Let $\sigma_p(A)$
(resp., $\sigma_a(A)$) denote the point spectrum (resp., the
approximate point spectrum) of the operator $A$, and let
$\rm{acc}\sigma(A)$ denote the set of  accumulation points of
$\sigma(A)$.
\begin{thm}\label{thm1} If, for an operator $A\in\b$, there exists a
compact operator $K\in\b$ such that $A+K\in\P$, then $A+K$ has SVEP
at points $\lambda\in\Phi_{sf}(A)$. Consequently,
$\Phi^+_{sf}(A)=\{\lambda\in\Phi_{sf}(A):\ind(A-\lambda)>
0\}=\emptyset$.\end{thm}
\begin{demo} Suppose to the contrary that $A+K$ does not have SVEP
at a point $\lambda\in\Phi_{sf}(A)=\Phi_{sf}(A+K)$. Recall from
\cite[Theorem 3.23]{A} that if an operator $T\in\b$ has SVEP at a
point $\mu\in\Phi_{sf}(T)$, then $\mu\in\i\sigma_a(T)$. Hence, since
$A+K$ does not have SVEP at $\lambda\in\Phi_{sf}(A+K)$, we must have
$\lambda\in\rm{acc}\sigma_p(A+K)$. Consequently, there exists a
sequence $\{\lambda_n\}\subset\sigma_p(A+K)$ of non-zero eigenvalues
of $A+K$ converging to $\lambda$. Choose $\alpha$,
$\beta\in\{\lambda_n\}$, and let  $M$ denote the subspace generated
by the eigenvectors $(A+K-\alpha)^{-1}(0)\cup
\beta(A+K-\beta)^{-1}(0)$. Then $A_1= (A+K)|_M$ is a polaroid
operator with $\sigma(A_1)=\{\alpha, \beta\}$, which implies that
$(A_1-\alpha)^{-1}(0)$ and $(A_1-\beta)^{-1}(0)$ are mutually
orthogonal spaces (in the sense of G. Birkhoff: A subspace $M$ of
$\X$ is orthogonal to a subspace $N$ of $\X$ if $||m||\leq ||m+n||$
for every $m\in M$ and $n\in N$ \cite[P. 93]{DS}). Now choose a
$\lambda_m\in\{\lambda_n\}$. Then the mutual orthogonality of the
eigenspaces corresponding to distinct (non -trivial) eigenvalues
implies $$\rm{dist} (x, (A+K-\lambda)^{-1}(0)) \geq 1$$ for every
unit vector $x\in (A+K-\lambda_m)^{-1}(0)$. Define
$\delta(\lambda_m, \lambda)$ by $$\delta(\lambda_m, \lambda)=
\rm{sup}\{\rm{dist} (x, (A+K-\lambda)^{-1}(0)): x\in
(A+K-\lambda_m)^{-1}(0), ||x||=1\}.$$ Then $\delta(\lambda_m,
\lambda)\geq 1$ for all $m$, and
$${\frac{|\lambda_m-\lambda|}{\delta(\lambda_m, \lambda)}}
\longrightarrow 0 \hspace{2mm}\mbox{as}\hspace{2mm} m\rightarrow
\infty,$$ i.e., the reduced minimum modulus function satisfies
$$\gamma(A+K-\lambda)= {\frac{|\lambda_m-\lambda|}{\delta(\lambda_m,
\lambda)}} \longrightarrow 0 \hspace{2mm}\mbox{as}\hspace{2mm}
m\rightarrow \infty.$$ Since this implies $(A+K-\lambda)(\X)$ is not
closed, we have a contradiction (of our assumption
$\lambda\in\Phi_{sf}(A+K)$). Hence $A+K$ has SVEP at every
$\lambda\in\Phi_{sf}(A+K)=\Phi_{sf}(A)$. The fact that
$\Phi^{+}_{sf}(A)=\emptyset$ is now a straightforward consequence of
``$A$ has SVEP at $\lambda\in\Phi_{sf}(A)$ implies
$\ind(A-\lambda)\leq 0$". \end{demo}

It is well known (indeed, easily proved) that $$A\in \b\cap\p
\Longleftrightarrow A^*\in B(\X^*)\cap\p.$$ This equivalence does
not extend to $\P$ operators. To see this, consider an operator
$A\in\b$ such that both $A$ and $A^*$ are in $\P$. Then both $A$ and
$A^*$ have SVEP at points in $\Phi_{sf}(A)$ ($=\Phi_{sf}(A^*)$) by
the preceding theorem. Hence, for every such operator $A$,
$$\lambda\in\Phi_{sf}(A)\Longrightarrow \lambda\in\Phi_w(A) =
\{ \lambda: \lambda\in\Phi(A),  \ind(A-\lambda)=0\}.$$ But then
$\lambda\in\Phi_{sf}(A)\cap\sigma(A)$ is (an isolated point of
$\sigma(A)$ which happens to be) a finite rank pole of $A$. That
this is (in general) false follows from a consideration of the
forward unilateral shift $U\in\B$ (which is trivially $\P$ and
satisfies $\lambda\in\Phi_{sf}(U)$ for all $|\lambda|< 1$).

 \

We consider next a sufficient condition for $A+K\in\P$ to imply
$A+K\in\P$. If $A\in\b$ and $M$ is an invariant (assumed, as before,
to be closed) subspace of $A$, then $A$ has an upper triangular
matrix representation  $$A= \left(\begin{array}{clcr}A_1 & *\\0 &
A_2\end{array}\right) \in B(M\oplus M^{\perp})$$ with main diagonal
$(A_1, A_2)$. Generally, $\sigma(A)\subseteq \sigma(A_1)\cup
\sigma(A_2)$ and $\sigma_w(A)\subseteq \sigma_w(A_1)\cup
\sigma_w(A_2)$: Indeed,
\begin{eqnarray*} & & \sigma(A_1)\cup \sigma(A_2)=
\sigma(A)\cup\{\sigma(A_1)\cap\sigma(A_2)\}\hspace{2mm}\mbox{and}\\
& & \sigma_w(A_1)\cup \sigma_w(A_2)=
\sigma_w(A)\cup\{\sigma_w(A_1)\cap\sigma_w(A_2)\}.\end{eqnarray*}
Recall from \cite[Exercise 7, P. 293]{TL} that
$$\asc(A_1-\lambda)\leq \asc(A-\lambda)\leq \asc(A_1-\lambda)+
\asc(A_2-\lambda)$$ for every complex $\lambda$; hence, if
$H_0(A-\lambda)= (A-\lambda)^{-p}(0)$ for a $\lambda\in
\{\sigma(A)\cap\sigma(A_1)\}$ (and some integer $p\geq 1$), then
\begin{eqnarray*}  H_0(A_1-\lambda) & = & H_0(A-\lambda)|_M \subseteq
(A-\lambda)^{-p}(0)\cap M= (A-\lambda)^{-p}(0)|_M\\ & = &
(A_1-\lambda)^{-p}(0)\subseteq H_0(A_1-\lambda);
\end{eqnarray*} consequently  $$H_0(A_1-\lambda)=
(A_1-\lambda)^{-p}(0)$$ (with $\asc(A_1-\lambda)\leq p$).

\

Let $\Pi_0(A)$ denote the set of Riesz points (i.e., finite rank
poles), and let $\Pi(A)$ denote the set of poles, of $A\in\b$. If
$A$ has SVEP on the complement of $\sigma_w(A)$ in $\sigma(A)$, then
$$ \sigma(A)\setminus\sigma_w(A)=\Pi_0(A) \Longleftrightarrow
\sigma(A)\setminus\sigma_{Bw}(A)=\Pi(A)$$ \cite[Theorem 2.1]{AZ}.

\begin{thm}\label{thm2} If, for an operator $A\in\b$, there exists a
compact operator $K\in\b$ such that $A+K\in\p$ and if
$\sigma(A+K)\setminus\sigma_w(A+K)=\Pi_0(A+K)$, then a sufficient
condition for $A+K\in\P$ is that:\\  (i)
$\sigma_w(A_1)\cup\sigma_w(A_2)\subseteq\sigma_w(A+K)$ for every
invariant subspace $M$ of $A+K$ such that $(A+K)|_M=A_1$ (and $(A_1,
A_2)$ is the main diagonal in the upper triangular representation of
$A+K\in B(M\oplus M^{\perp}))$;\\  (ii) $\i{\sigma_w(A_1)}\subseteq
\i{\sigma_w(A)}$.\end{thm} \begin{demo} We claim that
$\sigma(A+K)=\sigma(A_1)\cup\sigma(A_2)$. To prove the claim, we
start combining hypothesis (i) with the observation that
$\sigma_w(A+K)\subseteq\sigma_w(A_1)\cup\sigma_w(A_2)$ for every
upper triangular operator with main diagonal $(A_1, A_2)$ to obtain
$\sigma_w(A+K)=\sigma_w(A_1)\cup\sigma_w(A_2)$. Consider a complex
$\lambda\notin\sigma(A+K)$. Since \begin{eqnarray*} A+K-\lambda=
\left(\begin{array}{clcr} 1 & 0 \\0 & A_2-\lambda\end{array}\right)
\left(\begin{array}{clcr} 1 & * \\0 & 1\end{array}\right)
\left(\begin{array}{clcr} A_1-\lambda & 0 \\0 &
1\end{array}\right),\end{eqnarray*} $A_1-\lambda$ is left
invertible, $A_2-\lambda$ is right invertible,
$\alpha(A_1-\lambda)=0=\beta(A_2-\lambda)$ and $\ind(A+K-\lambda)=$
($\ind(A_1-\lambda)+\ind(A_2-\lambda)=0 \Longrightarrow$)
$\beta(A_1-\lambda)=\alpha(A_2-\lambda)$. If $\beta(A_1-\lambda)\neq
0$, then $\lambda\in\sigma_w(A_1)\cup\sigma_w(A_2)=\sigma_w(A+K)$, a
contradiction (since $\lambda\notin\sigma(A+K)$). Consequently,
$\beta(A_1-\lambda)=\alpha(A_2-\lambda)=0$, and hence
$\lambda\notin\{\sigma(A_1)\cup\sigma(A_2)$. This proves our claim.
Consider now a $\lambda\in\i\sigma(A_1)$. Then either
$\lambda\notin\sigma_{Bw}(A_1)$ or $\lambda\in\sigma_{Bw}(A_1)$. If
$\lambda\notin\sigma_{Bw}(A_1)$ and $\lambda\in\i\sigma(A_1)$, then
$\lambda$ is a pole of $A_1$. If, instead,
$\lambda\in\sigma_{Bw}(A_1)\subseteq\sigma_w(A_1)$, then
$\lambda\in\i{\sigma_w(A_1)}\subseteq\i{\sigma_w(A)}=\i{\sigma_w(A+K)}$.
Since $\sigma_w(A+K)=\sigma(A+K)\setminus{\Pi_0(A+K)}$,
$\lambda\in\i{\sigma_w(A+K)}$ implies $\lambda\in\i{\sigma(A+K)}$
and $\lambda\notin{\Pi_0(A+K)}$. Again, since $A+K\in\p$,
$\lambda\in{\Pi(A+K)}$ (is a pole of $A+K$ of infinite
multiplicity), and there exists an integer $p> 0$ such that
$H_0(A+K-\lambda)=(A+K-\lambda)^{-p}(0)$. But then, as seen above,
$H_0(A_1-\lambda)=(A_1-\lambda)^{-p}(0)$, and hence $\lambda$ is a
pole of $A_1$. This contradiction implies
$\lambda\notin\sigma_{Bw}(A_1)$, and $A_1\in\P$.\end{demo}

\begin{rema}\label{rema0} {\em The hypothesis
$\sigma_w(A_1)\cup\sigma_w(A_2)\subseteq\sigma_w(A+K)$ is not
necessary in Theorem \ref{thm2}. For example, if
$\sigma_w(A_1)\subseteq\sigma_w(A)$, then $\lambda\notin\sigma(A+K)$
implies $\alpha(A_1-\lambda)=0=\beta(A_2-\lambda)$ and
$\beta(A_1-\lambda)=\alpha(A_2-\lambda)$; hence, since
$\lambda\notin\sigma_w(A_1)$,
$\beta(A_1-\lambda)=\alpha(A_2-\lambda)=0$. Consequently,
$\sigma(A+K)=\sigma(A_1)\cup\sigma(A_2)$. The hypothesis
$\sigma(A+K)=\sigma(A_1)\cup\sigma(A_2)$ on its own does not
guarantee $A+K\in\P$ in Theorem \ref{thm2}. Let $R\in\b$ be a Riesz
operator and let $Q\in\b$ be a compact quasinilpotent operator.
Define $A\in B(\X\oplus\X)$ by $A=R\oplus 0$, and let $K=0\oplus Q$.
Then $A+K$ is a Riesz operator. Since the restriction of a Riesz
operator to an invariant subspace is again a Riesz operator
\cite{{H}, {A}}, $\sigma_w(A_1)\subseteq\sigma_w(A+K)$ for every
part (i.e., restriction an invariant subspace) $A_1=(A+K)|_M$ of
$A+K$. Hence $\sigma(A+K)=\sigma(A_1)\cup\sigma(A_2)$ for every
upper triangular representation, with main diagonal $(A_1, A_2)$, of
$A+K$. Evidently, $A+K\in\p$. Observe however that
$\i{\sigma_w(A_1)}\subseteq \i{\sigma_w(A+K)}$ fails for the (upper
triangular matrix) representation $A+K=Q\oplus A$ of $A+K$. Clearly,
$A+K\notin\P$. In the presence of the hypothesis
$\sigma(A+K)=\sigma(A_1)\cup\sigma(A_2)$, a sufficient condition for
$A_1\in\P$ is (of course) that
$\i\sigma(A_1)\cap\rm{acc}\sigma(A)=\emptyset$.}\end{rema}

{\bf Hilbert Space Operators.} Given an operator $A\in\B$, there
always exists a compact operator $K\in\B$ satisfying (the hypotheses
of Theorem \ref{thm2} that)
$\sigma(A+K)\setminus\sigma_w(A+K)=\Pi_0(A+K)$ and $A+K\in\p$: This
follows from the following familiar (see \cite{{He}, {ZL}, {LZ}})
argument. Every $A\in\B$ has an upper triangular matrix
representation
$$A=\left(\begin{array}{clcr}A_0 & *\\0 & A_1\end{array}\right)\in
 B(\H_0\oplus\H_1),\hspace{1mm} \H_1=\H\ominus\H_0,
\hspace{1mm}\sigma(A_0)=\Pi_0(A),\hspace{1mm}\sigma(A_1)=\sigma(A)\setminus{\Pi_0(A)}.$$
Consider $A_1\in B(\H_1)$. If we define $d$ by
$d=\rm{max}\{\rm{dist}(\lambda,\partial{\Phi_{sf}(A_1)}:\lambda\in\Pi_0(A)\}
< {\epsilon/2}$ (for some arbitrarily small $\epsilon> 0$), then
there exists a compact operator $K_1\in\ B(\H_1)$,
$||K_1||<{\epsilon/2}+d<\epsilon$, such that
$\rm{min.}\ind(A_1+K_1-\lambda)=0$ for all
$\lambda\in\Phi_{sf}(A_1)$ and $\sigma(A_1+K)=\sigma_w(A_1)$
\cite[Theorem 3.48]{He}. Let $A_{11}=A_1+K_1$. Then
$\lambda\in\i\sigma_w(A_1)$ and $\lambda\notin\sigma_{ulf}(A_{11})$
implies $\lambda\in\Pi_0(A)$; hence
$\i\sigma(A_{11})\cap\sigma_{ulf}(A_{11})\neq\emptyset$. Let
($\emptyset\neq$) $\Gamma\subset
\{\i\sigma(A_{11})\cap\sigma_{ulf}(A_{11})\}$. Then, for every
$\epsilon> 0$, there exists a compact operator $K_{11}\in B(\H_1)$,
$||K_{11}||< \epsilon$, such that
$$A_{11}+K_{11}=\left(\begin{array}{clcr}N & C\\0 &
A_2\end{array}\right)\in B(\H_{11}\oplus\H_{12}),
\hspace{2mm}\H_{11}=\H_1\ominus\H_{12}, \hspace{2mm}
\rm{dim}(\H_{11})=\infty,$$ $N$ is a diagonal normal operator of
uniform infinite multiplicity, $\sigma(N)=\sigma_{ulf}(N)=\Gamma$,
$\sigma(A_2)=\sigma(A_{11})$,
$\sigma_{ulf}(A_2)=\sigma_{ulf}(A_{11})$,
$\ind(A_2-\lambda)=\ind(A_{11}-\lambda)$ and
$\rm{min}.\ind(A_2-\lambda)=0$ for all $\lambda\in\Phi_{sf}(A_{11})$
\cite[Lemma 2.10]{Ji}. Assume, without loss of generality, that
$N=\oplus_{i=1}^{\infty}{\lambda_i 1_{\H_{11i}}}\in
B(\oplus_{i=1}^{\infty}{\H_{11i}})= B(\H_{11})$, where
$\rm{dim}(\H_{11i})=\infty$ for all $i\geq 1$. The points
$\lambda_i$ being isolated in $\sigma(N)$, there exists $\epsilon>
0$, an $\epsilon$-neighbourhood ${\cal N}_{\epsilon}(\lambda_i)$ of
$\lambda_i$ and a sequence $\{\lambda_{ij}\}\subset {\cal
N}_{\epsilon}(\lambda_i)$ such that $|\lambda_{ij}-\lambda_i|<
{\epsilon/2^i}$ for all $i\geq 1$. Choose an orthonormal basis
$\{e_{ij}\}_{j=1}^{\infty}$ of $\H_{11i}$, and let $K_{i0}$ be the
compact operator
$$K_{i0}=\sum_{j=1}^{\infty}{(\lambda_{ij}-\lambda_i)(e_{ij}\otimes e_{ij})}\in
B(H_{11i}), \hspace{1mm} ||K_{i0}||=
\rm{max}_j{|\lambda_{ij}-\lambda_i|},$$  define the compact operator
$K_{22}$ by$$K_{22}=\oplus_{i=1}^{\infty}{K_{i0}}\in B(\H_{11}),$$
and let $$N+K_{22}= \oplus_{i=1}^{\infty}{\{\lambda_i 1_{\H_{11i}}+
K_{i0}\}}= \oplus_{i=1}^{\infty}{N_i}.$$ Then each $N_i$ is a
diagonal operator with $\rm{diag}\{\lambda_{i1}, \lambda_{i2},
\cdots \}$, $\sigma(N+K_{22})=\cup_{i=1}^{\infty}{\sigma(N_i)}$ and
$\sigma_{ulf}(N+K_{22})$ is the closure of the set $\{\lambda_i:
1\leq 1\}$. Define the compact operator $K\in B(H_0\oplus \H_1)$ by
$$K=0\oplus(K_1+K_{11}+(K_{22}\oplus 0))\in B (\H_0\oplus
\H_{11}\oplus \H_{12}),$$ and consider a point
$\lambda\in\i\sigma(A+K)$: Either $\lambda\in\sigma(A_0)$, in which
case $\lambda\in{\Pi_0(A)}$, or,
$\lambda\in\i{\sigma_w(A+K)}=\i{\sigma_w(A)}$. If
$\lambda\in\i\sigma_w(A)$, then
$\lambda\in\i{\sigma_{ulf}(A)}=\i{(\Gamma)}$. Consequently,
$\lambda=\lambda_i$ for some integer $i\geq 1$, which then forces
$\lambda_i=\lim_{j\rightarrow\infty}{\lambda_{ij}}$. Since this
contradicts $\lambda\in\i\sigma(A+K)$, we are led to conclude
$A+K\in\p$.

\

The operator $\left(\begin{array}{clcr}A_0 & *\\0 &
N+K_2\end{array}\right)\in B (\H_0\oplus\H_{11})$ (of the above
construction) has SVEP. However, since
$\rm{min}.\ind(A_2-\lambda)=0$ for all $\lambda\in\Phi_{sf}(A_2)$,
either $\alpha(A_2-\lambda)=0$ or $\alpha(A_2-\lambda)^*=0$. If
$\alpha(A_2-\lambda)=0$, then ($\ind(A_2-\lambda)< 0$, and hence)
$A_2^*$ does not have SVEP at $\overline{\lambda}$; if, instead,
$\alpha(A_2-\lambda)^*=0$, then ($\ind(A_2-\lambda)> 0$, and hence)
$A_2$ does not have SVEP at $\lambda$. Conclusion: {\em The operator
$A+K$ above does not always satisfy the necessary condition of
Theorem \ref{thm1}, and hence  $A+K$ may or may not satisfy
$A+K\in\P$}. It is extremely complicated, if not impossible, to
determine the structure of the invaraint subspaces of the operator
$A+K$, and as such the determination of the passage from $A+K\in\p$
to $A+K\in\P$ does not seem to be within reach. An amenable case is
the one in which $A+K$ satisfies
$\Phi_{sf}^+(A+K)=\Phi_{sf}^+(A)=\emptyset$. Recall, \cite[Theorem
6.4]{He}, $A\in\B$ is {\em quasitriangular} if and only if
$\Phi_{sf}^+(A^*)=\emptyset$; if $A$ is quasitriangular, then there
is a compact operator $K$ such that $A+K$ is triangular. Thus, if
$\Phi_{sf}^+{A}=\emptyset$, then there exists a compact operator
$K\in\B$ and an orthonormal basis $\{e_i\}_{i=1}^{\infty}$ such that
$$(A+K)^*= \left(\begin{array}{clcr}a_{11} & a_{12} & \cdots\\0 &
a_{22} & \cdots\\\vdots & \vdots &
\ddots\end{array}\right)\left(\begin{array}{clcr}e_1\\e_2\\\vdots\end{array}\right),$$
for some scalars $a_{ii}\neq a_{jj}$ for all $i\neq j$. For each
invariant subspace $M$ of $A+K$ $$A+K=\left(\begin{array}{clcr}A_1 &
*\\0 &
A_2\end{array}\right)\left(\begin{array}{clcr}M\\M^{\perp}\end{array}\right)
\Longleftrightarrow (A+K)^*=\left(\begin{array}{clcr}A_2^* & *\\0 &
A_1^*\end{array}\right)\left(\begin{array}{clcr}M^{\perp}\\M\end{array}\right),$$
where $A_1^*$ has an upper triangular matrix with main diagonal
$\rm{diag}(A_1^*)=\{a_{{n_k}{n_k}}\}_{k=1}^{\infty}$,
$\sigma(A+K)^*=\sigma(A_1^*)\cup\sigma(A_2^*)$ and
$\iso\sigma(A_1^*)\subset\i\sigma(A+K)^*$. Applying Theorem
\ref{thm2} (to obtain \cite[Theorem 5.1]{LZ}) and combining with
Theorem \ref{thm1} we have: \begin{thm}\label{thm3} Given an
operator $A\in\B$, a necessary and sufficient condition that there
exist a compact operator $K\in\B$ such that $A+K\in\P$ is that
either (i) $\Phi_{sf}^+(A)=\emptyset$ or (equivalently) (ii) $A+K$
has SVEP at points in $\Phi_{sf}(A)$.\end{thm}

 \

 We close this note with the result that given an operator $A\in\B$, a necessary and sufficient condition
 for there to exist a compact operator $K\in\B$ such that $A+K\in\P$ is that
 the component $\Omega_a(A)$
 is connected. Here, given an operator $A$, the
 component $\Omega_a(A)$ of $\Phi_{sf}(A)$  is defined by
 $$ \Omega_a(A)=\{\lambda\in\Phi_{sf}(A):\ind(A-\lambda)\le 0\}.$$  \begin{thm}\label{thm4}  (i) If,
 for an operator $A\in\b$, the component $\Omega_a(A)$ of
 $\Phi_{sf}(A)$ is connected, then $A+K$ has SVEP on $\Phi_{sf}(A)$
 for every compact operator $K\in\b$.\\ \noindent (ii) If $\X=\H$ is a
 Hilbert space and $A\in\B$, then a necessary and sufficient condition
 for there to exist a compact operator $K\in\B$ such that $A+K\in\P$ is that
  the component $\Omega_a(A)$ of $\Phi_{sf}(A)$
 is connected.\end{thm} \begin{demo}  {\em (i)}  We prove by contradiction. If $\Omega_a(A)$ is connected, then (it
 has no bounded component, and hence) it has just one component,
 namely itself, and hence the resolvent set $\rho(A)$ intersects $\Omega_a(A)$. Consequently,
 both $A$ and $A^*$ have SVEP at points in $\Omega(A)$ \cite[Theorem
 3.36]{A}. Suppose now that there exists a compact operator $K\in\b$
 such that $A+K$ does not have SVEP at a point
 $\lambda\in\Phi_{sf}(A+K)=\Phi_{sf}(A)$. Since $(A+K)^*$ has SVEP
 and $A+K$ fails to have SVEP at a point  $\lambda\in\Phi_{sf}(A)$ implies
 $\ind(A-\lambda)>0$, we must have that neither of $A+K$ and
 $(A+K)^*$ have SVEP at $\lambda$. Hence
 $\asc(A+K-\lambda)=\dsc(A+K-\lambda)=\infty$. On the other hand,
 since $\rho(A+K)\subset\Omega_a(A)$, the continuity of the index at
 points $\lambda\in\Omega_a(A)$ implies that $\ind(A+K-\lambda)=0$.
 Thus $\alpha(A+K-\lambda)=0$ (except perhaps for a countable set of
 $\lambda$), and it follows that $A+K-\lambda$ is bounded below
 (and hence $\asc(A+K-\lambda)<\infty$).This is a
 contradiction.

\

\noindent{\em (ii)} Start by observing that if $A+K$ has SVEP at
$\lambda\in\Phi_{sf}(A)$, then (necessarily) $\ind(A+K-\lambda)\leq
0$, equivalently, $\Phi_{sf}^+(A)=\emptyset$, for every compact
operator $K\in\B$. This, by Theorem \ref{thm3} above or
\cite[Theorem 5.1]{LZ}, implies the existence of a compact operator
$K\in\B$ such that $A+K\in\P$. Conversely, if there exists a compact
operator $K$ such that $A+K\in\P$, then $A+K$ has SVEP on
$\Phi_{sf}(A)$. Assume, to the contrary, that $\Omega_a(A)$ (is not
connected, and hence) has a bounded component $\Omega_0(A)$. Then
$\Gamma=\partial{\Omega_0}(A)\subset\sigma_{ulf}(A)$, and there
exists a compact operator $K_1\in\B$ such that $A+K_1$ has the upper
triangular matrix representation $$ A+K_1=\left(\begin{array}{clcr}
N & *\\0 & A_2\end{array}\right)\in B(\H_1\oplus \H_2),
\hspace{2mm}\rm{dim}(\H_1)=\infty,$$ with respect to some
decomposition $\H=\H_1\oplus\H_2$ of $\H$, where $N$ is a normal
diagonal operator of uniform infinite multiplicity,
$\sigma(N)=\sigma_{uf}(N)=\Gamma$, $\sigma_{uf}(A_2)=\sigma_{uf}(A)$
and $\ind(A_2-\lambda)=\ind(A-\lambda)$ for all
$\lambda\in\Phi_{sf}(A)$ (see \cite[Lemma 2.10]{Ji}). The spectrum
$\sigma(N)=\Gamma$ of $N$ being the boundary of a bounded connected
open subset of $\C$, \cite[Theorem 3.1]{He1} implies the existence
of a compact operator $K_2\in B(\H_1)$ such that $\sigma(N+K_2)$
equals the closure $\overline{\Omega_0(A)}$. Define the compact
operator $K\in B(\H_1\oplus\H_2)$ by $K=K_1+(K_2\oplus 0)$. Then
$$A+K=\left(\begin{array}{clcr}N+K_2 & *\\0 & A_2\end{array}\right)\in
\B,$$ where for every $\mu\in\Omega_0(A)$ we have
$\mu\in\Phi_{sf}(N+K_2)=\Phi_{sf}(N)$ with $\ind(N+K_2-\mu)=0$. It
being clear that SVEP for $A+K$ at a point implies SVEP for $N+K_2$
at the point, it follows that every $\mu\in\Omega_0(A)$ is an
isolated point -- a contradiction. Conclusion: $\Omega_a(A)$, has no
bounded component.  \end{demo}

\vskip10pt \noindent\normalsize\rm B.P. Duggal, 8 Redwood Grove,
Northfield Avenue, Ealing, London W5 4SZ, United Kingdom.\\
\noindent\normalsize \tt e-mail: bpduggal@yahoo.co.uk

\end{document}